\documentclass{amsart}

\usepackage{macros}
\standardmargins\standardpackages\standardrenews\standardlabeling\standardmakes
\colorcommentstrue

\draftfalse

\authormumtaz
\authordavid

\begin{document}
\title{A general principle for Hausdorff measure}

\begin{abstract}
We introduce a general principle for studying the Hausdorff measure of limsup sets. A consequence of this principle is the well-known Mass Transference Principle of Beresnevich and Velani (2006). 
\end{abstract}
\maketitle

                \section{Introduction and Statements of main results}

A fundamental problem in the theory of metric Diophantine approximation is to determine the `size' of the limsup set $$\limsup_{i\to\infty}B_i=\bigcap_{n=1}^\infty\bigcup_{i=n}^\infty B_i=\{x\in X: x\in B_i \text{ for infinitely many } i\in \mathbb N\}$$  in terms of Lebesgue measure, Hausdorff dimension, or Hausdorff measure. Here and throughout  we let  $(B_i)_i$ be a sequence of open sets in a metric space $X$. Let $f$ be a dimension function i.e. an increasing continuous function 
$f:[0, \infty) \to [0, \infty)$ with $f(0)=0$. We denote by $\mathcal H^f$ the $f$-dimensional Hausdorff measure, which is proportional to the standard Lebesgue measure when $X = \mathbb R^d$ and $f(r)=r^d$. In the case where the dimension function is of the form $f(r):=r^s$ for some $s > 0$, $\mathcal H^f$ is simply denoted as $\mathcal H^s$. For the definitions of Hausdorff measure and dimension and their properties we refer to the book \cite{Falconer_book2013}. The Hausdorff--Cantelli lemma \cite[Lemma 3.10]{BernikDodson} states that  $\mathcal H^f(\limsup_{i\to\infty} B_i)=0$ if $\sum_{i=1}^\infty f({\diam}(B_i))<\infty$, where $\diam(B_i)$ denotes the diameter of $B_i$. Here the emphasis is on using a particular `nice' cover of the limsup set, namely $\{B_i:i\in\mathbb N\}$, to establish an upper bound for Hausdorff measure.  In contrast, proving the $f$-dimensional Hausdorff measure to be positive is a challenging task, requiring all possible coverings to be considered  and, therefore, represents the main problem in metric Diophantine approximation (in various settings).  

\begin{question}\label{q1} Under what conditions is $\mathcal H^f(\limsup_{i\to\infty}B_i)$ strictly positive?

\end{question}

The following principle commonly known as the Mass Distribution Principle  \cite[\64.1]{Falconer_book2013}  has been the go-to method in giving an answer to Question \ref{q1}. 

\begin{lemma}[Mass Distribution Principle]\label{mdp} Let $\mu$ be a probability measure supported on a subset $F$ of $X$. Suppose there are positive constants $c>0$ and $\epsilon>0$ such that $$\mu (U)\leq c f(\diam(U))$$ for all sets $U$ with $\diam (U)\leq \epsilon$. Then $\mathcal H^f(F)\geq \mu(F)/c.$

\end{lemma}

Specifically, the mass distribution principle replaces the consideration of all coverings by the construction of a particular measure $\mu$. Given a sequence of sets $(B_i)_i$, if we want to prove that $\mathcal H^f(\limsup_{i\to\infty} B_i)$ is infinite (and in particular strictly positive), one possible strategy is to deploy the mass distribution principle in two steps:

\begin{itemize}
\item construct a suitable Cantor type subset $\mathcal K\subset F = \limsup_{i\to\infty} B_i$ and a probability measure $\mu$ supported on $\mathcal K$,

\item  show that for any fixed $c>0$, $\mu$  satisfies the condition that for any measurable set $U$ of sufficiently small diameter,  $\mu (U)\leq c f(\diam(U))$. 
 \end{itemize}
If this can be done, then by the mass distribution principle, it follows that  $$\mathcal H^f(F)\geq \mathcal H^f(\mathcal K)\geq c^{-1}.$$ Then since $c$ is arbitrary, it follows that $\mathcal H^f(F)=\infty$.

The main intricate and substantive part  of this entire process is the construction of a suitable Cantor type subset of $F$ which supports a probability measure. We introduce a generalised principle to determine the $f$-dimensional Hausdorff measure of limsup sets which throws out the Cantor type construction from this process.

To state our result, we introduce some notation. Let $X$ be a metric space.   For $\delta>0$, a measure $\mu$ is {\em Ahlfors $\delta$-regular} if and only if there exist positive constants $0<c_1<1< c_2<\infty$ and $r_0 > 0$ such that the inequality 
$$c_1r^\delta\leq \mu (B(x,r)) \leq c_2 r^\delta$$ 
holds for every ball $B:=B(x, r)$ in $X$ of radius $r \leq r_0$ centred at $x \in \Supp(\mu)$, where $\Supp(\mu)$ denotes the topological support of $\mu$. The space $X$ is called Ahlfors $\delta$-regular if there is an Ahlfors $\delta$-regular measure whose support is equal to $X$. If $X$ is Ahlfors $\delta$-regular, then so is the $\mathcal H^\delta$ measure restricted to $X$ i.e.  $\mathcal H^\delta\given X$.  We will frequently be using the notation $B^f:=B(x, f(r)^{1/\delta})$. For real quantities $A,B$ that depend on parameters, we write $A \lesssim B$ if $A \leq c B$ for a constant $c > 0$ that is independent of those parameters. We write  $A\asymp B$ if $A\lesssim B\lesssim A$.


\begin{theorem}\label{mainthm}
Fix $\delta > 0$, let $(B_i)_i$ be a sequence of open sets in an Ahlfors $\delta$-regular metric space $X$, and let $f$ be a dimension function such that
\begin{align} \label{f1}
r \mapsto &r^{-\delta} f(r) \text{ is decreasing, and}\\ \label{f2}
&r^{-\delta} f(r) \to \infty \text{ as }r \to 0.
\end{align}
Fix $C > 0$, and suppose that the following hypothesis holds:
\begin{itemize}
\item[(*)] For every ball $B_0 \subset X$ and for every $N\in\mathbb N$, there exists a probability measure $\mu = \mu(B_0,N)$ with $\Supp(\mu) \subset \bigcup_{i\geq N} B_i\cap B_0$, such that for every ball $B = B(x,\rho) \subset X$, we have
\begin{equation}
\label{muBbound}
\mu(B) \lesssim \max\left(\left(\frac{\rho}{\diam B_0}\right)^\delta,\frac{f(\rho)}{C}\right).
\end{equation}
\end{itemize}
Then for every ball $B_0$,
\[
\mathcal H^f\left(B_0 \cap \limsup_{i\to\infty} B_i\right) \gtrsim C.
\]
In particular, if the hypothesis \text{(*)} holds for all $C$, then
\[
\mathcal H^f\left(B_0\cap \limsup_{i\to\infty} B_i\right) = \infty.
\]
\end{theorem}

The condition \eqref{f2} is a natural condition which implies that $\mathcal H^f(B)=\infty$. The hypothesis (*) is the  main ingredient of this theorem and, roughly speaking, this gives a systematic way of constructing the probability measure on the limsup set. In contrast,  however, the mass distribution principle (Lemma 1.2) assumes that there is a mass distributed on the limsup set. 

The hypothesis essentially translates to saying that for any given ball $B_0$ with a probability measure $\mu$, we can find a bunch of smaller balls $B_i$'s which are inside the approximation $\bigcup_{i\geq N} B_i\cap B_0$ to the limsup set $\limsup_{i\to \infty} B_i$ with a mass assigned to each of them. Basically we need to assume that we can, for each ball $B_0$, redistribute the mass on smaller balls $B_i$ in a way that whenever we draw a test ball $B$ then the measure of this ball can always be bounded by either the number of full balls inside it  i.e. $f(\rho)/C$ or  by the measure of the fractions of balls $\left(\frac{\rho}{\diam B_0}\right)^\delta$.

\subsection{The Mass Transference Principle}\label{MTP}

In a landmark paper \cite{BeresnevichVelani}, Beresnevich and Velani  introduced the Mass Transference Principle which has become a major tool in converting Lebesgue measure theoretic statements for limsup sets into Hausdorff measure statements for limsup sets. This is surprising as the Lebesgue measure is the `coarser' notion of `size' than the Hausdorff measure.  As one would expect, the Mass Transference Principle has many applications in number theory, such as derivation of Jarn\'ik's theorem from Khintchine's theorem or the derivation of the Jarn\'ik--Besicovitch theorem from Dirichlet's theorem. Other than that the Mass Transference Principle can be used to determine Hausdorff measure and dimension of limsup sets in the context of dynamical systems such as $\beta$-dynamical systems \cite{CHW, LuWu}. This principle has further been generalised to 
lim sup sets defined via neighbourhoods of sets satisfying a certain local
scaling property satisfying the open set
condition in \cite{AllenBaker} and to the multifractal formalisms in \cite{FST}.

\begin{theorem}[Beresnevich--Velani, 2006]\label{mtp} Let $X\subset \mathbb R^d$ be Ahlfors $\delta$-regular. Let $(B_i)_{i\in\mathbb N}$ be a sequence of balls in $X$ with $\textrm{rad}(B_i)\to 0$ as $i\to \infty$. Let $f$ be a dimension function such that $r\mapsto r^{-\delta}f(r)$ is monotonic. Suppose that for every ball $B\subset X$

\begin{equation}\label{eq1mtp}\mathcal H^\delta(B\cap\limsup_{i\to\infty}B_i^f)=\mathcal H^\delta(B).\end{equation}
Then for every ball $B\subset X$
\begin{equation}\label{eq2mtp}\mathcal H^f(B\cap\limsup_{i\to\infty}B_i)=\mathcal H^f(B).\end{equation}
\end{theorem}

The Mass Transference Principle has been significantly used in determining the Hausdorff measure, and as a consequence the Hausdorff dimension, of limsup sets. For instance, using the Mass Transference Principle the Hausdorff measure version of the Duffin-Schaeffer conjecture was formulated in \cite{BeresnevichVelani}. By considering $X$ to be the middle third Cantor set, for which $\delta=\log 2/\log 3$, the Mass Transference Principle has been used in \cite{LSV} in establishing a complete metric theory for sets of $\psi$-approximable points and as a consequence proving the existence of very well approximable numbers other than Liouville numbers in the middle third Cantor set, an assertion that was attributed to Mahler.

\medskip

\subsection{Theorem \ref{mainthm} $\Longrightarrow$ Theorem \ref{mtp}}

Before establishing this implication we state the following covering lemma, which is a variant of the $K_{G,B}$ lemma of Beresnevich and Velani \cite{BeresnevichVelani}. The difference is that we replace the set $B_i^f$ by the set $B_i^{f/C}$, where $C > 0$ is an arbitrary constant, while we need bounds that are independent of $C$. The proof is more involved, and is a variant of the proof of the Vitali covering lemma \cite[Theorem 2.2]{Mattila}.

\begin{lemma}\label{kgblemma}
Fix $C > 0$. Let $(B_i = B(x_i,\rho_i))_{i\in\mathbb N}$ be a sequence of balls in an Ahlfors $\delta$-regular space $X$ with $\textrm{rad}(B_i)\to 0$ as $i\to \infty$. Let $f$ be a dimension function such that \eqref{eq1mtp} is satisfied for every ball $B$ in $X$. Fix $B_0 \subset X$ and $N\in\mathbb N$. Then there exists a finite set $I = I(B_0,N) \subset \{N,\ldots\}$ such that the collection $\left(B_i^{f/C}\right)_{i\in I}$ is disjoint and its union is a subset of $B_0$ satisfying
\begin{equation}
\label{KGB}
\mathcal H^\delta\left(\bigcup_{i\in I} B_i^{f/C} \right) \geq \frac12 \mathcal H^\delta(B_0).
\end{equation}
\end{lemma}
Note that it is necessary in the later argument to have a factor of $\frac12$ in \eqref{KGB}, rather than the factor of $\frac1C$ that would result from an easier proof.
\begin{proof}
Without loss of generality assume that $N = 1$. By rearranging we can without loss of generality assume that $\rho_j \leq 2\rho_i$ whenever $j \geq i$. Let $I^\infty \subset \mathbb N$ be defined using a greedy algorithm: $j\in I^\infty$ if and only if $B_i^{f/C} \cap B_j^{f/C} = \emptyset$ for all $i < j$ such that $i\in I^\infty$. Fix $M$ and let
\[
U = \mathrm{interior}(B_0) \butnot \bigcup_{\substack{i\leq M \\ i\in I^\infty}} B_i^{f/C}.
\]
Fix $x\in U\cap \limsup_{i\to\infty} B_i^f$. Then there exists sufficiently large  $i > M$ such that $x\in B_i^f \subset U$. Since $U$ is an open set which means that there is a ball around $x$ contained in $U$ which in return implies that for sufficiently large $i$, $B_i^{f/C}\subset U$.  If $i\notin I^\infty$, then there exists $j < i$ such that $j\in I^\infty$ and $B_i^{f/C}\cap B_j^{f/C} \neq \emptyset$, otherwise, simply take $j=i$. In other words, there exists $j\leq i$ such that $j\in I^\infty$ and $\rho_i\leq 2\rho_j$ and $B_i^{f/C}\cap B_j^{f/C}\neq \emptyset$.  

Since $B_i^{f/C}\subset U$, it follows that $j > M$. It follows from \ref{f1} that $f$ is doubling i.e. there exists a constant $\lambda>1$ such that for $\rho>0$, $f(2 \rho)\leq \lambda^\delta f(\rho)$. Now since $\rho_i \leq 2\rho_j$, $f$ is doubling and from the well-known Vitali's $5r$-covering lemma, we have that 

$$B_i^f \subset \eta B_j^f,$$
where $\eta=\lambda^\delta 5$.

Thus,
\[
U\cap \limsup_{i\to\infty} B_i^f \subset \bigcup_{\substack{i > M \\ i\in I^\infty}} \eta B_i^f.
\]
Since \eqref{eq1mtp} is satisfied, by Ahlfors regularity we have
\[
\mathcal H^\delta(U) \leq \sum_{\substack{i > M \\ i\in I^\infty}} \mathcal H^\delta(\eta B_i^f) \asymp_C \sum_{\substack{i > M \\ i\in I^\infty}} \mathcal H^\delta(B_i^{f/C}) = \mathcal H^\delta\left(U \cap \bigcup_{i\in I^\infty} B_i^{f/C}\right).
\]
Since $M$ was arbitrary, it follows that
\[
\mathcal H^\delta\left(B_0 \cap \bigcup_{i\in I^\infty} B_i^{f/C}\right) = \mathcal H^\delta(B_0)
\]
and thus by continuity of measures, taking an appropriate finite initial segment of $I^\infty$ completes the proof.
\end{proof}

Let $I = I(B_0,N)$ be as in Lemma \ref{kgblemma}, and construct the measure $\mu = \mu(B_0,N)$ as follows:
\[
\mu = \frac1K \sum_{i\in I} \mathcal H^\delta\left(B_i^{f/C}\right) \frac{\mathcal H^\delta\given B_i}{\mathcal H^\delta(B_i)}
\]
where $K$ is a constant chosen so that $\mu(X) = 1$, i.e. $K = \sum_{i\in I} \mathcal H^\delta(B_i^{f/C})$. Note that by \eqref{KGB},
\[
K \asymp \mathcal H^\delta(B_0).
\]
Also note that by \eqref{f2}, by increasing $N$ if necessary we can assume without loss of generality that $f(\rho_i)/C \geq (2\rho_i)^\delta$, and thus that $B_i \subset B_i^{f/C}$, for all $i\geq N$. Since $B_i^{f/C} \subset B_0$ for all $i\in I$, we have $$\Supp(\mu) \subset \bigcup_{i\in I} B_i \subset \bigcup_{i\geq N} B_i\cap B_0.$$ Fix $B = B(x,\rho)$; to complete the proof, we need to show that \eqref{muBbound} holds. Suppose first that $B$ intersects only one element of the collection $(B_i)_{i\in I}$, say $B_i = B(x_i,\rho_i)$. Then
\begin{align*}
\mu(B) &= \frac1K \mathcal H^\delta\left(B_i^{f/C}\right) \frac{\mathcal H^\delta(B_i\cap B)}{\mathcal H^\delta(B_i)}\\
&\lesssim \frac1{\mathcal H^\delta(B_0)} \frac{f(\rho_i)}{C} \min\left(1,\left(\frac{\rho}{\rho_i}\right)^\delta\right)
\by{Ahlfors $\delta$-regularity}\\
&\leq \frac1{\mathcal H^\delta(B_0)} \frac{f(\rho)}{C}. \by{\eqref{f1}}
\end{align*}
On the other hand, suppose that $B$ intersects multiple elements of $(B_i)_{i\in I}$. Since $\left(B_i^{f/C}\right)_{i\in I}$ is a disjoint collection, it follows that for all $i\in I$ for which $B_i\cap B \neq \emptyset$, we have $$2\rho \geq \diam(B) \geq \dist\left(B_i,X\butnot B_i^{f/C}\right) \geq (1/2)(f(\rho_i)/C)^{1/\delta}.$$ It follows that $\diam(B_i^{f/C}) \leq 8\rho$ and thus
\[
B_i^{f/C} \subset B(x,9\rho)
\]
so
\begin{align*}
\mu(B) &\leq \frac1K \sum_{\substack{i\in I \\ B\cap B_i \neq \smallemptyset}} \mathcal H^\delta\left(B_i^{f/C}\right)\\
&\leq \frac1K \mathcal H^\delta(B(x,9\rho)) \since{$(B_i^{f/C})_{i\in I}$ are disjoint} \\
&\asymp \left(\frac{\rho}{\diam(B_0)}\right)^\delta. \note{Ahlfors $\delta$-regularity}
\end{align*}
This completes the proof.

\section {Proof of Theorem \ref{mainthm}}

For each $n\in\mathbb N$, we construct a set $T^n \subset \mathbb N^n$, a family of balls $(B_\omega)_{\omega\in T^n}$, and a family of parameters $(p_\omega)_{\omega\in T^n}$, recursively as follows:
\begin{itemize}
\item For $n=0$, we let $T^0 = \{\emptyset\}$, $B_\smallemptyset = B(x_\smallemptyset,1)$, and $p_\smallemptyset = 1$. Here $\emptyset\in \mathbb N^0$ is the empty string, and $x_\smallemptyset\in X$ is chosen arbitrarily.
\item Fix $n\in\mathbb N$, and suppose that $T^n$ and the families $(B_\omega)_{\omega\in T^n}$ and $(p_\omega)_{\omega\in T^n}$ have been defined. Fix $\omega\in T^n$, and let $\mu_\omega = \mu(B_\omega,n)$, where the notation $\mu(B,N)$ is as in hypothesis (*). Choose
\begin{equation}
\label{rhoomegadef}
0 < \rho_\omega < \dist\Big(\Supp(\mu_\omega),X\butnot \bigcup_{i\geq n} B_i\cap B_\omega\Big) \text{ such that} \left(\frac{\rho_\omega}{\diam B_\omega}\right)^\delta \leq \frac{f(\rho_\omega)}{C}.
\end{equation}
Such a choice is possible by \eqref{f2}.

Let $F_\omega \subset K_\omega = \Supp(\mu_\omega)$ be a maximal $4\rho_\omega$-separated set (Here, a set $F$ is called \emph{$\rho$-separated} if for all $u_1, u_2 \in F$, we have $d(u_1, u_2)\geq \rho$).  Let $\pi_\omega: K_\omega\to F_\omega$ be a map such that $\dist(x,\pi_\omega(x)) \leq 4\rho_\omega$ for all $x\in K_\omega$. Let $k_\omega = \#(F_\omega)$, and let $(x_{\omega a})_{1 \leq a \leq k_\omega}$ be an enumeration of $F_\omega$. Finally, let
\[
T^{n+1} = \{ \omega a : \omega \in T^n, \; 1 \leq a \leq k_\omega\}.
\]
For each $\omega a \in T^{n+1}$, let
\[
B_{\omega a} = B(x_{\omega a},\rho_\omega) \;\;\;\; \text{ and } \;\;\;\;
p_{\omega a} = p_\omega \cdot (\pi_\omega)_*[\mu_\omega](x_{\omega a}),
\]
where  $(\pi_\omega)_*[\mu_\omega]$ denotes the pushforward of a measure $\mu_\omega$ under the map $\pi_\omega$.
\end{itemize}

Note that for all $\omega\in T = \bigcup_{n\geq 0} T^n$, we have
\[
\sum_{a \leq k_\omega} p_{\omega a} = p_\omega
\]
and for all $a\leq k_\omega$ we have
\[
B_{\omega a} \subset \bigcup_{i\geq |\omega|} B_i \cap B_\omega
\]
where $|\omega|$ denotes the length of $\omega$. Moreover, by the definition of $F_\omega$, the collection $(B_{\omega a})_{a\leq k_\omega}$ is disjoint. It follows that there is a unique probability measure $\mu$ such that
\[
\mu(B_\omega) = p_\omega \text{ for all } \omega\in T,
\]
and that $\Supp(\mu) \subset \limsup_{i\to \infty} B_i$. To complete the proof, we need to show that
\begin{equation}
\label{massdist}
\mu(B(x,\rho)) \lesssim f(\rho)/C
\end{equation}
for any sufficiently small ball $B(x,\rho)$. This will allow us to apply the mass distribution principle (Lemma \ref{mdp}).

Indeed, fix such a ball $B = B(x,\rho)$, let $\omega\in T$ be the shortest word such that $B\cap B_\omega \neq \emptyset$ and $\rho_\omega \leq \rho$, and let $\tau$ be the initial segment of $\omega$ that is one letter shorter than $\omega$, so that $\rho_\tau > \rho$. Since $F_\tau$ is $4\rho_\tau$-separated, for all $\tau a\neq \omega$ we have $\dist(x_\omega,x_{\tau a}) \geq 4\rho_\tau$ and thus $$\dist(B_\omega,B_{\tau a}) \geq 2\rho_\tau > 2\rho,$$ from which it follows that $B\cap B_{\tau a} = \emptyset$. So $\mu(B) = \mu(B\cap B_\omega)$, and thus
\begin{align*}
\mu(B) &\leq \sum_{\substack{a \leq k_\omega \\ B\cap B_{\omega a} \neq \smallemptyset}} p_{\omega a}\\
&\leq p_\omega \sum_{\substack{a \leq k_\omega \\ x_{\omega a} \in B(x,2\rho)}} (\pi_\omega)_*[\mu_\omega](x_{\omega a})\\
&\leq p_\omega \mu_\omega\left(\bigcup_{\substack{a \leq k_\omega \\ x_{\omega a}\in B(x,2\rho)}} B(x_{\omega a},4\rho_\omega)\right)\\
&\leq p_\omega \mu_\omega\big(B(x,6\rho)\big)\\
&\lesssim p_\omega \max\left(\left(\frac{\rho}{\diam B_\omega}\right)^\delta,\frac{f(\rho)}{C}\right)
\end{align*}
by our assumption on $\mu_\omega$. A similar argument gives
\[
p_\omega = \mu(B_\omega) \lesssim p_\tau \max\left(\left(\frac{\rho_\tau}{\diam B_\tau}\right)^\delta,\frac{f(\rho_\tau)}{C}\right).
\]
Combining with \eqref{rhoomegadef} and using the fact that $p_\tau \leq 1$ gives
\[
p_\omega \lesssim \frac{f(\rho_\tau)}{C}
\]
and thus since $p_\omega \leq 1$,
\[
\mu(B) \lesssim \max\left(\left(\frac{\rho}{\diam B_\omega}\right)^\delta \frac{f(\rho_\tau)}{C},\frac{f(\rho)}{C}\right).
\]
Applying \eqref{f1} with $r = \rho_\tau$ and $ar = \rho$ and using the fact that $\diam B_\omega = 2\rho_\tau$ demonstrates \eqref{massdist} and completes the proof.

\medskip

\noindent{\bf Acknowledgments.} The first-named author was supported by La Trobe University's start-up grant.  The second-named author was supported by the EPSRC Programme Grant EP/J018260/1. We would like to thank Demi Allen for several conversations regarding the contents of this paper and to the anonymous referee for useful comments.

\providecommand{\bysame}{\leavevmode\hbox to3em{\hrulefill}\thinspace}
\providecommand{\MR}{\relax\ifhmode\unskip\space\fi MR }
\providecommand{\MRhref}[2]{%
  \href{http://www.ams.org/mathscinet-getitem?mr=#1}{#2}
}
\providecommand{\href}[2]{#2}

\end{document}